%% file: Ops_and_Hypothesis_K_starMay_ArXiV.tex
\documentclass{amsart}
\usepackage{amsmath, amssymb}

\begin{document}
	\newcommand{\xtra}[1]{}
	\newcommand{\xtrasec}[1]{}
	\newcommand{\ques}[1]{{\bfseries \emph{#1}}}

\input{format}

\title[Discrete Fractional Integral Operators]{On Discrete Fractional Integral Operators and \\
Mean Values of Weyl Sums}
\author{Lillian B. Pierce}
\address{School of Mathematics, Institute for Advanced Study, Princeton New Jersey 08540}
\email{lbpierce@math.ias.edu}
\thanks{The author was supported by the Simonyi Fund at the Institute for Advanced Study and the National Science Foundation, including DMS-0902658 and DMS-0635607, during this research.}

\date{21 May 2010}
\maketitle

\begin{abstract}
In this paper we prove new $\ell^p \maps \ell^q$ bounds for a discrete fractional integral operator by applying techniques motivated by the circle method of Hardy and Littlewood to the Fourier multiplier of the operator.  
From a different perspective, we describe explicit interactions between the Fourier multiplier and mean values of Weyl sums. These mean values express the average behaviour of the number $r_{s,k}(l)$ of representations of a positive integer $l$ as a sum of $s$ positive $k$-th powers. Recent deep results within the context of Waring's problem and Weyl sums enable us to prove a further range of complementary results for the discrete operator under consideration. 
\end{abstract}

\section{Introduction}

Let $I_{k,\lam}$ be the operator acting on (compactly supported) functions $f: \Z \maps \C$ by
\beq\label{Ik_dfn}
 I_{k,\lam} f(n) = \sum_{m=1}^\infty \frac{f(n-m^k)}{m^\lam},
 \eeq
where $0< \lam < 1$ and $k \geq 1$ is an integer. For which pairs of exponents $p,q$ may $I_{k,\lam}$ be extended to a bounded operator from $\ell^p(\Z)$ to $\ell^q(\Z)$? Here  we denote by $\ell^p(\Z)$ the space of functions $f: \Z \maps \C$ such that the norm $||f||_{\ell^p(\Z)} = (\sum_{n\in \Z} |f(n)|^p)^{1/p}$ is finite. 

One may formulate the expected mapping properties of $I_{k,\lam}$  by 
considering its continuous analogue, defined for $0<\lam <1$ by
\[ \mathcal{I}_{k,\lam}f(x)=\int_1^\infty f(x-y^k) y^{-\lam}dy = \frac{1}{k}\int_{1}^\infty f(x-u)u^{(1-\lam)/k -1}du.\]
By the classical Hardy-Littlewood-Sobolev fractional integration theorem, $\mathcal{I}_{k,\lam}$ is known
to be bounded from $L^p(\R)$ to $L^q(\R)$ if and only if $1< p< q < \infty$ and $1/q = 1/p - (1-\lam)/k$.
In analogy, one would expect the following for the discrete operator $I_{k,\lam}$:
\begin{conjecture}\label{Conj}
For $0< \lam < 1$ and $k \geq 1$, $I_{k,\lam}$ extends to a bounded operator from $\ell^p(\Z)$ to $\ell^q(\Z)$ if and only if $p,q$ satisfy
\newline
 (i) $1/q \leq 1/p - (1-\lam)/k$,
\newline
 (ii)  $1/q<\lam, 1/p >1-\lam.$ 
\end{conjecture}

That these two conditions are necessary for $I_{k,\lam}$ to be bounded from $\ell^p$ to $\ell^q$ may be seen by considering two simple examples: for condition (i), set $f(n)=n^{-\gamma}$ if $n\geq 1$ and $f(n)=0$ if $n \leq 0$, where $\ga >1/p$: for condition (ii), set $g(0)=1$ and $g(n)=0$ if $n \neq 0$.

For $k=1$, the fact that conditions (i) and (ii) are also sufficient for $I_{k,\lam}$ to be bounded from $\ell^p$ to $\ell^q$ is an immediate consequence of the known $(L^p, L^q)$ bounds for the operator $\Ical_{k,\lam}$ in the continuous setting (see Proposition \emph{a} of \cite{SW2}). For any $k \geq 2$, a similar comparison to the continuous operator $\Ical_{k,\lam}$ shows that $I_{k,\lam}$ maps $\ell^p$ to $\ell^q$ for $1/q \leq 1/p -(1-\lam)$ (see Proposition \emph{b} of \cite{SW2}). However, note that for $k \geq 2$, this is far from sharp, compared to condition (i) of Conjecture \ref{Conj}.

Stein and Wainger \cite{SW2} \cite{SW3}, Oberlin \cite{Obe}, and Ionescu and Wainger \cite{IW} have studied the discrete operator $I_{k,\lam}$  extensively in the case $k=2$, ultimately proving the full sharp bounds of Conjecture \ref{Conj} in this specific case. (Indeed, in \cite{IW} Ionescu and Wainger further prove a deep result for translation invariant discrete singular Radon transforms that implies that for $\lam=1+i\ga$ with $\ga \neq 0$, $I_{k,\lam}$ is a bounded operator on $\ell^p(\Z)$ for all $1<p<\infty$.)
For $k \geq 3,$ Conjecture \ref{Conj} remains open, and the partial results obtained so far have been slight. 
In this paper we prove new bounds toward Conjecture \ref{Conj} for $I_{k,\lam}$ for all $k \geq 3$, and along the way, we outline precisely why the higher degree problem appears to be difficult.

We will take two distinct approaches to the problem. Define 
the discrete Fourier transform of a function $f \in \ell^1(\Z )$ by
\[ \hat{f}(\theta) = \sum_{n \in \Z} f(n) e^{-2\pi i n \theta}.\]
The operator $I_{k,\lam}$ is translation invariant, and as a result $(I_{k,\lam}f)\hat{\;}(\theta) = m_{k,\lam}(\theta) \hat{f}(\theta)$, where a simple computation shows that the Fourier multiplier $m_{k,\lam}$ is defined by
\[ m_{k,\lam}(\theta) = \sum_{n=1}^\infty \frac{e^{-2\pi i n^k \theta}}{n^\lam},\]
for $\theta \in [0,1]$.

We will observe that to prove Conjecture \ref{Conj}, it would be sufficient to show that for $0< \lam<1$, $m_{k,\lam}$ belongs to the weak-type $L^r$[0,1]  space\footnote{The weak-type $L^r$ space is the space of functions such that $|\{\theta : |f(\theta)| >\al \}| \leq c \al^{-r}$ for all $\al >0$, which when normed becomes the Lorentz space $L^{r,\infty}$.} for $r =k/(1-\lam)$. In the first half of this paper we will prove a result of this type by means of a decomposition of the multiplier $m_{k,\lam}$ motivated by the circle method of Hardy and Littlewood. Estimates for exponential sums play a key role in bounding $m_{k,\lam}$, and we apply bounds due to refinements of the Weyl bound and Vinogradov's mean value theorem. The results we obtain via this method (Theorem \ref{thm_circle}) prove the results of Conjecture \ref{Conj} for a restricted range of $\lam$ (approximately of size $1- \frac{1}{\frac{3}{2} k^2 \log k} < \lam < 1$).

In fact, this limitation on $\lam$ reflects a profound difficulty embedded in this approach. 
For any $1/2 < \lam < 1$, proving that $m_{k,\lam} \in L^{r,\infty}[0,1]$ with $r=k/(1-\lam)$ is stronger than showing $m_{k,\lam} \in L^{2k}[0,1]$.\xtra{In fact, $L^r[0,1] \subset L^{r,\infty}[0,1] \subset L^p[0,1]$ if $r>p$.
For in fact, if $m \in L^{r,\infty}[0,1]$, set $\Lambda(\al) = |\{x: |m(x)| > \al\}|$, so that by assumption $\Lambda(\al) \leq A\al^{-r}$ for every $\al>0$. Then $||m||_{L^p[0,1]} =\frac{1}{p}\int_0^\infty \al^{p-1}\Lambda(\al)d\al = \int_0^1 + \int_1^\infty.$ Since the underlying space $[0,1]$ has finite measure, trivially $\Lambda(\al) \leq 1$, and hence the first integral is bounded. By the weak-type assumption on $\Lambda(\al)$ the second integral is bounded by $\int_1^\infty \al^{p-r-1}d\al,$ which is finite since $r>p$.} 
But as was first observed in \cite{SW2}, the statement that $m_{k,\lam}\in L^{2k}[0,1]$ for all $1/2 < \lam< 1$ is equivalent to Hypothesis $K^*$, a classical conjecture about the average behaviour of $r_{k,k}(l)$, the number of representations of a positive integer $l$ as a sum of $k$ $k$th positive powers. This conjecture remains unproved for $k \geq 3$. In the second half of this paper we examine a more general equivalence between properties of the multiplier $m_{k,\lam}$ and the average behaviour of $r_{s,k}(l)$, the number of representations of a positive integer $l$ as a sum of $s$ $k$th positive powers. For $s$ sufficiently small or sufficiently large with respect to $k$, best possible results for the relevant averages of $r_{s,k}(l)$ are known due to recent deep results on mean values of Weyl sums \cite{Ford} \cite{SalWoo}. These results will ultimately allow us to prove a collection of new bounds for the operator $I_{k,\lam}$ (Theorems \ref{I_sk_small}, \ref{I_kk}, \ref{I_sk}).

\subsection{Summary of results}
We note that  these two main approaches---$L^{r,\infty}[0,1]$ bounds for the Fourier multiplier $m_{k,\lam}$ coming directly from the circle method, versus $L^{2k}$ bounds coming from mean values of exponential sums and asymptotics for Waring's problem (and hence also ultimately from the circle method)---are not redundant, but complementary. Because the key distinctions between the four theorems we prove will be most clearly understood within the context of their proofs, we defer stating the precise theorems until they arise naturally in the course of the discussion. Instead, we provide a brief summary of the results, for large $k$. Theorems \ref{thm_circle} -- \ref{I_sk} prove that $I_{k,\lam}$ is a bounded operator from $\ell^p(\Z)$ to $\ell^q(\Z)$ under the following conditions, respectively:
\begin{itemize}
\item[(1)] $1/q \leq 1/p - \frac{1-\lam}{k}$ and $1- \frac{1}{\frac{3}{2} k^2 \log k}<\lam <1$;
\item[(2)]  $1/q < 1/p - \frac{4\log \log k(1-\lam)}{\log k}$ and $1/2 < \lam < 1$;
\item[(3)]  $1/q < q/p - \frac{(1-\lam)}{\log_2 k +1}$ and $ 1- \frac{\log_2 k +1}{2k}< \lam < 1$;
\item[(4)] $1/q< 1/p - \frac{1-\lam}{k}$ and $1 - \frac{1}{2k \log k} < \lam < 1$.
\end{itemize}

These results may be compared most intuitively by regarding the pictorial representation of $(\ell^p,\ell^q)$ bounds for the operator $I_{k,\lam}$ for each fixed $0<\lam<1$ as a set of points $(1/p, 1/q)$ in the unit square $[0,1] \times [0,1]$. In this interpretation, Conjecture 1 states that for each $0<\lam<1$, $I_{k,\lam}$ is expected to be bounded from $\ell^p$ to $\ell^q$ for all points $(1/p,1/q)$ that lie on or below the diagonal line $1/q = 1/p - (1-\lam)/k$ and also lie inside the open box described by $1/q<\lam,$ $1/p>1-\lam$. Generally speaking, the first approach (Theorem \ref{thm_circle}) gives $(\ell^p, \ell^q)$ bounds for $(1/p,1/q)$ lying on or below the diagonal line $1/q = 1/p - (1-\lam)/k$ but for a short range of $\lam$, while the second approach (Theorems \ref{I_sk_small}, \ref{I_kk}, \ref{I_sk}) proves $(\ell^p, \ell^q)$ bounds for $(1/p,1/q)$ further restricted to lying strictly below a lower diagonal line, but for a longer range of $\lam$. A pictorial representation summarising the distinct conditions for $p,q,\lam$ arising in each of Theorems \ref{thm_circle} to \ref{I_sk} is given in Figure 1.

\setlength{\unitlength}{1cm}

\begin{figure}\label{fig_1}
\begin{picture}(12,5)(0,0)

\thicklines
\put(0,0){\line(1,0){4}}
\put(0,0){\line(0,1){4}}
\put(4.25,-.1){$1/p$}
\put(-.2,4.25){$1/q$}


\put(1.2,0){\line(1,1){2.8}}
\thinlines
\put(1.3,0){\line(1,1){2.7}}
\put(1.9,0){\line(1,1){2.1}}
\put(2.2,0){\line(1,1){1.8}}
\put(4.2,2.8){(1), (4)}
\put(4.2,2.2){(3)}
\put(4.2,1.85){(2)}

\put(2.4,0){\line(0,1){1.6}}
\put(2.4,1.6){\line(1,0){1.6}}
\put(2.1,-.3){$1-\lam$}
\put(4.2,1.5){$\lam$}

\put(7,1){\line(1,0){4}}
\put(7,1.5){\line(1,0){4}}
\put(7,2){\line(1,0){4}}
\put(7,2.5){\line(1,0){4}}

\thicklines
\put(9,1){\line(1,0){2}}
	\put(9,1.02){\line(1,0){2}}
\put(10.3,1.5){\line(1,0){.7}}
	\put(10.3,1.52){\line(1,0){.7}}
\put(10.5,2){\line(1,0){.5}}
	\put(10.5,2.02){\line(1,0){.5}}
\put(10.8,2.5){\line(1,0){.2}}
	\put(10.8,2.52){\line(1,0){.2}}

\thinlines
\put(9,.95){\line(0,1){.12}}
\put(10.3,1.45){\line(0,1){.12}}
\put(10.5,1.95){\line(0,1){.12}}
\put(10.8,2.45){\line(0,1){.12}}

\put(7,.95){\line(0,1){.12}}
\put(7,1.45){\line(0,1){.12}}
\put(7,1.95){\line(0,1){.12}}
\put(7,2.45){\line(0,1){.12}}

\put(11,.95){\line(0,1){.12}}
\put(11,1.45){\line(0,1){.12}}
\put(11,1.95){\line(0,1){.12}}
\put(11,2.45){\line(0,1){.12}}

\put(11.5,.95){(2)}
\put(11.5,1.45){(3)}
\put(11.5,1.95){(4)}
\put(11.5,2.45){(1)}

\put(6.9,.5){0}
\put(10.9,.5){1}
\put(8.5,.5){$0<\lam<1$}

\end{picture}
\vspace{.2cm}
\caption[Figure 1]{This is a pictorial representation of the results of Theorems \ref{thm_circle} -- \ref{I_sk} (labelled by (1) -- (4)). The left plot represents those $p,q$ for which each theorem shows that $I_{k,\lam}$ is a bounded operator from $\ell^p(\Z)$ to $\ell^q(\Z)$. The right plot depicts a comparison of the ranges of $\lam$ in which each theorem holds.
}
\end{figure}
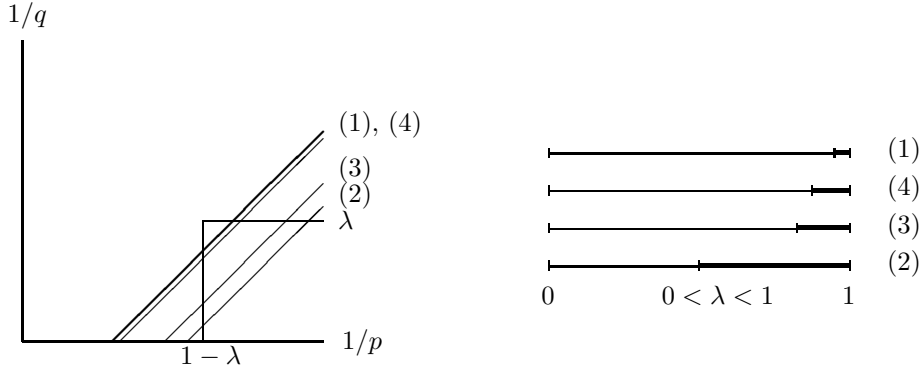

\subsection{Outline of the paper}
We will proceed as follows. In Section \ref{sec_prelim} we prove two technical lemmas relating properties of the multiplier $m_{k,\lam}$ to bounds for the operator $I_{k,\lam}$. In Section \ref{sec_weak_bds} we prove a weak-type bound for $m_{k,\lam}$ via a circle method decomposition of the multiplier, and indicate the resulting bound for the operator $I_{k,\lam}$. In Section \ref{sec_Hyp_K} we prove an explicit equivalence between $m_{k,\lam}$ and an average bound for $r_{s,k}(l)$, which we denote by Property $K_{s,k}^*(\be)$. Finally, in Section \ref{sec_mean_vals} we record the best known results for Property $K_{s,k}^*(\be)$ for three ranges of $s$ with respect to $k$, and their consequences for $I_{k,\lam}$. The author thanks Elias M. Stein for his interest and encouragement, Per Salberger and Trevor D. Wooley for providing a preprint of \cite{SalWoo}, and Tim Browning, Kevin Ford, and Roger Heath-Brown for helpful comments.

\section{Preliminaries}\label{sec_prelim}
We first record two functional analytic lemmas that will be used to derive results toward Conjecture \ref{Conj} from the knowledge that $m_{k,\lam} \in L^{r,\infty}[0,1]$ for $r = k/(1-\lam)$. The first result is a ``folk lemma'' noted by Stein and Wainger \cite{SW2}:

\begin{lemma}\label{folk_lemma}
Let $T$ be a convolution operator acting on functions $f: \Z \maps \C$ with Fourier multiplier $m$,
\[ (Tf)\hat{\;} (\theta) = m(\theta) \hat{f}(\theta).\]
If $m (\theta) \in L^{r,\infty}[0,1],$ then $T: \ell^p(\Z) \maps \ell^q(\Z) $ for $1/q = 1/p - 1/r$ and $1 < p \leq 2 \leq q < \infty$. 
\end{lemma}
For completeness, we record the proof of this lemma due to Stein and Wainger. First, assume that $q=2$, so that $1/p = 1/r + 1/2$. Then for $f \in \ell^p$, by Paley's version of the Hausdorff-Young inequality, $\hat{f} \in L^{p',p}[0,1]$, where $1/p + 1/p'=1$. Therefore by the multiplicative property of Lorentz spaces (see Chapter V \S 3 of \cite{SteinWeiss}),
\[ (Tf)\hat{\;}(\theta) = \hat{f}(\theta) m(\theta) \in L^{p',p} \cdot L^{r,\infty} \subset L^{p_0,q_0},\]
where $1/p_0 = 1/p' + 1/r = 1/2$, $1/q_0 = 1/p + 1/\infty = 1/p$. Therefore $(Tf)\hat{\;} \in L^{2,p} \subset L^{2,2} = L^2$ since $p \leq 2$ by assumption. Hence $Tf \in \ell^2$ and so $T$ maps $\ell^p$ to $\ell^2$. The case with $p=2$ and $1/2 - 1/r = 1/q$ follows by considering the adjoint operator to $T$, and the lemma then follows by interpolation between the two resulting bounds for $T$. 

Supposing that it is known that $m_{k,\lam} \in L^{r,\infty}[0,1]$ for $r=k/(1-\lam)$ in a range $\lam_k < \lam<1$, then Conjecture \ref{Conj} in this restricted range of $\lam$ follows immediately from Lemma \ref{folk_lemma}. For indeed, the restriction that $p\leq 2 \leq q$ can be removed by complex interpolation with the result of Ionescu and Wainger \cite{IW} that $I_{k,\lam}$ extends to a bounded operator on $\ell^p(\Z)$ for all $1< p< \infty$ on the line $\Re(\lam)=1$, $\Im(\lam) \neq 0$. The remaining $p,q$ pairs satisfying condition (i) in Conjecture \ref{Conj} with a strict inequality follow simply from the inclusion property $\ell^{q_1} \subset \ell^{q_2}$ if $q_1 < q_2$.

In fact, we note that it is not necessary to prove $m_{k,\lam} \in L^{r,\infty}[0,1]$ for all $0<\lam<1$ in order to prove Conjecture \ref{Conj} in full. It would suffice to show that $m_{k,\lam} \in L^{r,\infty}[0,1]$ for $\lam$ in the smaller range $\lam_k^*  < \lam < 1$, where $\lam_k^* = 1-\frac{k}{2k-1}$ is the crossover value at which condition (i) and condition (ii) (with equality) meet at a single point. For $\lam \leq \lam_k^*$, condition (ii) is stronger than condition (i). Thus for $\lam \leq \lam_k^*$, the desired $(\ell^p, \ell^q)$ bounds for $I_{k,\lam}$ could be obtained by interpolating the known result for $\Re(\lam) > \lam_k^*$ with the trivial $\ell^1 \maps \ell^\infty$ bound for $I_{k,\lam}$ for $\Re(\lam) \geq 0$, followed by applying the inclusion property of $\ell^q$ spaces, and taking adjoints. Asymptotically, $\lam_k^*$ approaches $1/2$ from below as $k \maps \infty$, so that the region of most interest in the work that follows is $1/2 < \lam < 1$.

The results we will prove in Section \ref{sec_mean_vals} will place $m_{k,\lam}$ not in a weak $L^r$ space varying with $\lam$, but in a fixed $L^p$ space; to translate these results into bounds for the operator $I_{k,\lam}$ we will apply the following lemma. 
\begin{lemma}\label{K_to_op}
Suppose $m_{k,\lam} \in L^{u}[0,1]$ for all $\eta <\lam < 1$ with $\eta >0$ and $u >2$ fixed. Then 
$ I_{k,\lam}$ extends to a bounded operator from $\ell^p(\Z)$ to  $\ell^q(\Z)$
if $\eta < \lam < 1$ and $p,q$ satisfy the conditions
\newline
(i) $1/q < 1/p - \frac{(1-\lam)}{u(1 - \eta)}$,
\newline
 (ii) $1/q <\lam, 1/p >1-\lam$.
  \end{lemma}
Indeed, since trivially $L^{u}[0,1] \subset L^{u,\infty}[0,1]$, under the hypotheses of the lemma it follows from Lemma \ref{folk_lemma} that $I_{k,\lam}$ maps $\ell^p$ to $\ell^q$ boundedly for $1/q = 1/p - 1/u$ for all $\eta<\lam<1$. (As before, we may remove the restriction $1 < p \leq 2 \leq q <\infty$ by interpolation with the $\ell^p \maps \ell^p$ result of \cite{IW} for $I_{k,\lam}$ on the line $\Re(\lam)=1$.) However, in this case note that $u$ is fixed (independent of $\lam$), so that this result is strongest for $\lam$ close to $\eta$. Thus the lemma follows from interpolating the $\ell^p \maps \ell^q$ result for $\lam$ arbitrarily close to $\Re(\lam)=\eta$ with the $\ell^p \maps \ell^p$ result on the line $\Re(\lam)=1$.\xtra{Note that because we may only take $\lam> \eta$, and not $\lam = \eta$, in Lemma \ref{mod_folk_lemma}, we do not obtain equality in condition (i) of the theorem.}

\section{Weak-type bounds for $m_{k,\lam}$ via the circle method}\label{sec_weak_bds}
Our first approach is to bound $m_{k,\lam}(\theta)$ directly by decomposing it according to the Diophantine properties of $\theta$. Our argument follows closely the work of Stein and Wainger \cite{SW2}; we improve on the original presentation by optimising the choice of the major and minor arcs and applying sharper bounds for Weyl sums that arise in the minor arcs. 

Decompose the multiplier as
\beq\label{j_decomp}
 m_{k,\lam}(\theta) = \sum_{n=1}^\infty \frac{e^{-2\pi i n^k \theta}}{n^\lam} 
	= \sum_{j=0}^\infty \sum_{2^j \leq n < 2^{j+1}} \frac{e^{-2\pi i n^k \theta}}{n^\lam}.
	\eeq
We define the major and minor arcs in terms of real parameters $\be, \be_0, \be_1 \geq0$ that will be fixed later. Given $j \geq 0$, Dirichlet's approximation principle guarantees that for each $\theta \in [0,1]$ there exist integers $a,q$ such that $1 \leq q \leq 2^{(\be - \be_1)j}, 1 \leq a \leq q$, $(a,q)=1$, with
\beq\label{Dir_theta}
 | \theta - a/q| \leq \frac{1}{q2^{(\be - \be_1)j}}.
 \eeq
For each $1 \leq q \leq \frac{1}{10}2^{\be_0 j}$ and $1 \leq a \leq q$ with $(a,q)=1$, we define the major arc 
\[ \M_j(a/q)  = \{ \theta : | \theta - a/q| \leq 1/(q2^{(\be - \be_1)j}) \}.\]
The minor arcs (for each fixed $j$) are then defined to be the complement of the union of the major arcs in $[0,1]$.

The key property of the major arcs is that they are disjoint if $(a_1,q_1) \neq (a_2,q_2)$. For indeed, if they were not we would have
\[ \frac{1}{q_1q_2} \leq \left| \frac{a_1}{q_1} - \frac{a_2}{q_2} \right| \leq \frac{1}{q_1 2^{(\be-\be_1)j}} +  \frac{1}{q_2 2^{(\be-\be_1)j}} ,\]
which is impossible for $1 \leq q_1, q_2 \leq \frac{1}{10}2^{\be_0j}$, as long as $0 \leq \be_0 \leq \be - \be_1$, which we henceforward assume.

 It will also be convenient to define a version of the major arcs that is independent of $j$; for every $(a,q)=1$ let
\[ \M^*(a/q)  =  \{ \theta : | \theta - a/q| \leq 1/(10q^2) \}.\]
Then certainly $\M_j(a/q) \subset \M^*(a/q)$ for every $j$, since $ \be_0 \leq \be - \be_1$. 
These intervals are also disjoint, as long as $q_1 \leq q_2 \leq 2q_1$, for assuming $\M^*(a_1/q_1)$ intersects $\M^*(a_2/q_2)$ with $(a_1,q_1) \neq (a_2,q_2)$, we would have
\[ \frac{1}{q_1q_2} \leq \left| \frac{a_1}{q_1} - \frac{a_2}{q_2} \right| \leq \frac{1}{10q_1^2} + \frac{1}{10q_2^2},\]
which is impossible under the assumed conditions.

\subsection{The major arcs}
We consider the inner sum on the right hand side of (\ref{j_decomp}). For $\theta \in \M_j(a/q)$, write $\theta=a/q+\al$ and write $n=mq+l$ in terms of its residue class modulo $q$. 
In the case $k=2$, $m_{k,\lam}$ is related to a theta function, and applying the Jacobi inversion formula (i.e. Poisson summation) separates out the arithmetic information and main term from this inner sum, with a small remainder term. For higher degrees, only cruder estimates are available.
First note that
	\beq\label{inner_dyad}
	\sum_{2^j \leq n < 2^{j+1}} \frac{e^{-2\pi i n^k \theta}}{n^\lam}
		= \sum_{l=1}^q e^{-2\pi i l^k a/q} \sum_{2^j/q \leq m < 2^{j+1}/q}  \frac{e^{-2\pi i (mq+l)^k \al}}{(mq+l)^\lam}+O(q2^{-j\lam}).\eeq
We apply the van der Corput summation formula (see Lemma 8.8 of \cite{IK}):
	\[ \sum_{m=a}^b e^{2\pi i f(m)}g(m) = \int_a^b e^{2\pi i f(x)} g(x) dx + O(2^{-j\lam}), \]
with $f(x) = (xq+l)^{k}\al$, $g(x) = (xq+l)^{-\lam}$ and $a=2^j/q$, $b=2^{j+1}/q$. 
Note that $f'(x)$ is monotonic and bounded by
\[ |f'(x)| \leq k|\al|(2^{j+1} +q)^{k-1} \leq k|\al|(2^{j+1} +2^{\be_0j})^{k-1} =O(2^{-(\be - \be_1)j} 2^{j(k-1)}),\]
as long as $\be_0 \leq 1$; this is uniformly bounded if $\be - \be_1 \geq k-1$, which we now assume.
	Thus it  suffices to consider the integral
	\[ \int_{2^j/q}^{2^{j+1}/q} e^{2\pi i f(x)} g(x) dx
		= \frac{1}{q}\int_{2^j + l}^{2^{j+1} +l} e^{2\pi i y^k \al} \frac{dy}{y^\lam} 
		= \frac{1}{q}\int_{2^j}^{2^{j+1}} e^{2\pi i y^k \al} \frac{dy}{y^\lam} +O(2^{-j\lam}).\]
This last expression is independent of $l$, which was our goal.

For $(a,q)=1$, define the Weyl sum
\[ S(a/q) = \sum_{l=1}^q e^{2\pi i l^k a/q},\]
which admits the classical bound $S(a/q) = O(q^{1-1/k+\ep})$ for any $\ep >0$ (see Lemma 3 in Chapter 2 of \cite{Vin}).  
We also define the integral function
\[ \Phi(u) = \int_{1}^2 e^{2\pi i y^k u} \frac{dy}{y^\lam}.\]
Then for $\theta \in \M_j(a/q)$, we have shown that
\[ \sum_{2^j \leq n < 2^{j+1}} \frac{e^{2\pi i n^k\theta}}{n^\lam} =
	 \frac{1}{q} S(a/q) 2^{j(1-\lam)} \Phi(2^{jk}\al) + O(q2^{-j\lam}).\]
Summing in $j$, the error term contributes
\[ O(\sum_j q2^{-j\lam}) = O(\sum_j 2^{(\be_0 -\lam)j}) = O(1),\]
assuming $ \lam > \be_0$. 
Next, note that 
\[ \sum_{j=0}^\infty 2^{j(1-\lam)} |\Phi(2^{jk}u)| = O(|u|^{-1/r}),\] 
where $r=k/(1-\lam)$; this follows from splitting the sum for each $u$ into those $j$ such that $2^{jk}|u| \leq 1$ and $2^{jk}|u| > 1$, and using the bounds $|\Phi(u)| \leq A$ and $|\Phi(u)| \leq A/|u|,$ respectively, for some constant $A$.

Thus in total, the contribution of the major arcs to $m_{k,\lam}$ is (up to a constant)
\beq\label{s_sum}
 \sum_{a,q} q^{-1/k+\ep}|\theta-a/q|^{-1/r}\chi_{a/q}(\theta) = \sum_{s=0}^\infty \sum_{2^s \leq q < 2^{s+1}} \sum_{\bstack{1 \leq a \leq q}{(a,q)=1}}	q^{-1/k+\ep}|\theta-a/q|^{-1/r}\chi_{a/q}(\theta) ,
 \eeq
where $\chi_{a/q}$ denotes the characteristic function of the interval $\M^*(a/q)$. Note that for each fixed $a,q$ pair, the function $q^{-1/k+\ep}|\theta-a/q|^{-1/r}\chi_{a/q}(\theta) $ has $L^{r,\infty}[0,1]$ norm $O(2^{-s(1/k-\ep)})$.  In order to sum for each fixed $s$ the $O(2^{2s})$ such functions in (\ref{s_sum}) and retain a finite $L^{r,\infty}$ norm, we apply the following lemma (Lemma 1 of \cite{SW2}):
\begin{lemma}
Given $N$ functions $f_1, \ldots, f_N$ with disjoint supports and with $f_j$ uniformly in $L^{r,\infty}[0,1]$, then the sum
\[ F_N = N^{-1/r} \sum_{j=1}^N f_j\]
belongs to $L^{r,\infty}[0,1]$ uniformly in $N$.
\end{lemma}\label{SW_lemma}
To see this, note that by assumption, for each $j$, $|\{ |f_j| >\al \}| \leq \al^{-r}$ for all $\al>0$. Thus by the disjoint support hypothesis,
\[ \{x: |F(x)| > \al \}  = \Union_{j=1}^N  \{x: N^{-1/r} |f_j(x)| > \al \},\]
and thus
\[|\{x: |F(x)| > \al \} | = \sum_{j=1}^N | \{x: N^{-1/r} |f_j(x)| > \al \} | \leq N^{-1} \sum_{j=1}^N \al^{-r}  = \al^{-r}.\]

As a result of Lemma \ref{SW_lemma}, (\ref{s_sum}) is bounded by
$\sum_s 2^{-s(1/k-\ep)} 2^{2s/r},$ where $r=k/(1-\lam)$,
which is finite if and only if $\lam >1/2$. 
To summarize, we have shown that the contribution of the major arcs to $m_{k,\lam}$ is in $L^{r,\infty}[0,1]$ under the assumptions
\beq\label{maj_assp}
 \be_0 \leq 1, \quad \be - \be_1 \geq k-1, \quad \lam > \be_0, \quad \lam>1/2.
 \eeq

\subsection{The minor arcs}
We now turn to the minor arcs, for which we must once again bound the inner sum in (\ref{j_decomp}), namely
\beq\label{inner}
\sum_{2^j \leq n < 2^{j+1}} \frac{e^{-2\pi i n^k \theta}}{n^\lam}.
\eeq
Let 
\[ S_{N}(\theta) = \sum_{M \leq n <M+ N} e^{-2\pi i n^k \theta}.\]
It is sufficient to bound $S_{N}(\theta)$ with $N=M=2^j$, 
for then by partial summation (\ref{inner}) is dominated by $O(2^{-j\lam}|S_{2^j}(\theta)|)$, and the total contribution of the minor arcs is bounded by
\beq\label{minor_total}
\sum_{j =0}^\infty 2^{-j\lam}\sup_{\theta \in \chi_{\mathfrak{m}_j}} |S_{2^j}(\theta)|,
\eeq
where $\chi_{\mathfrak{m}_j}$ is the characteristic function of the minor arcs for $j$.

In general, the Weyl bound and Vinogradov's mean value theorem provide bounds of the form $S_{N}(\theta) = O(N^{1-\sig_k})$ for some $\sig_k >0$ depending on $k$ and the Diophantine properties of $\theta$. 
In order to state the best known bounds for $S_N(\theta)$, we require notation relating to the multi-dimensional mean value of Weyl sums defined by
\beq\label{Jcal_dfn}
{\bf J}_{s,k}(X) = \int_{[0,1]^k} \left| \sum_{x=1}^X e(\al_1 x + \cdots \al_k x^k) \right|^{2s}  d\al_1 \cdots d\al_k.
\eeq
Simple estimates (see \cite{Ford}) show that ${\bf J}_{s,k}(X) \gg X^{2s-k(k+1)/2}$;
Vinogradov's mean-value methods prove upper bounds of the form
\beq\label{J_eta}
{\bf J}_{s,k}(X)   \ll X^{2s-k(k+1)/2 + \eta(s,k)},
\eeq
where $\eta(s,k)$ is small if $s$ is large enough relative to $k$. 
The current best known value of $\eta(s,k)$ is 
\beq\label{eta_sk}
\eta(s,k) \approx k^2 e^{-2s/k^2},
\eeq
 due to Wooley \cite{Woo92}, which is close to zero as soon as $s$ is of the order $k^2 \log k$. (For more precise statements of the best known upper bounds for $\eta(s,k)$, see Lemmas 5.1 and 5.5 of \cite{Ford}.)

The subject of extracting a bound for an individual exponential sum $S_N(\theta)$ from the average ${\bf J}_{s,k}(X)$ requires care; we record below the current best known bounds (also compiled in Lemma 5.3 of \cite{Ford}), under the assumption $\be_0=\be_1$.
\begin{lemma}\label{sig_lemma}
Suppose that there exist integers $a,q$ such that $|\theta-a/q| \leq 1/q^2$ where $1 \leq a \leq q$, $(a,q)=1$ and $N^{\be_0} < q \leq N^{k-\be_0}$. Then 
\[ |S_N(\theta) | \leq C_{k,\ep} N^{1-\sig_k(\be_0)+\ep} \]
where 
\[ \sig_{k}(\be_0) = \max \left(\sig_k^{(1)}, \sig_k^{(2)}, \sig_k^{(3)}\right)      \]
and 
\begin{eqnarray}
\sig_k^{(1)}& = &  \frac{\be_0}{2^{k-1}}, \label{prime1}\\
 \sig_k^{(2)} &=& \max_{s \geq 1} \left(\frac{\be_0 - \eta(s,k-1)}{2s}\right), \label{prime2} \\
 \sig_k^{(3)} & = & \max_{1 \leq r \leq k/2} (\min(\sig_k^{(3a)}(r), \sig_k^{(3b)}(r)) \label{prime3}
 \end{eqnarray}
 where under the assumption $\be_0 > 1-1/k$,
 \begin{eqnarray*}
 \sig_k^{(3a)}(r) & = & \max_{s \geq k(k-1)/2} \left( \frac{r - \eta(s,k-1)}{2rs} \right), \\
 \sig_k^{(3b)}(r) & = & \max_{t \geq 1} \left( \frac{k-r(1+\eta(t,k))}{2tk} \right).
 \end{eqnarray*}
\end{lemma}
The value $\sig_k^{(1)}$ is simply Weyl's bound (see for example Lemma 2.4 of Vaughan \cite{Vau81}). The value $\sig_k^{(2)}$ follows directly from a version of Vinogradov's mean value theorem (Theorem 5.3 of Vaughan \cite{Vau81}). The value $\sig_k^{(3)}$ also comes from Vinogradov's mean value theorem, in a refinement due to Theorem 2 of Wooley \cite{Woo95}; it is valid for all $\be_0 > 1-1/k.$

To make the results of Lemma \ref{sig_lemma} easier to interpret, we note that roughly speaking, with $\eta(s,k) \approx  k^2 e^{-2s/k^2}$ as in (\ref{eta_sk}), the optimal choice for $s$ in (\ref{prime2}) is approximately $s \approx 2 (k-1)^2 \log(k-1)$, leading to 
\[\sig_k^{(2)}\approx \frac{\be_0 -(k-1)^{-2}}{4(k-1)^2\log(k-1)}.\]
For large $k$, $\sig_k^{(3)}$ is the most powerful choice. In this direction,
by Corollary 1 of \cite{Woo95}, for any $\ep>0$ there exists $k(\ep)$ such that if $k\geq k(\ep)$ and $\be_0 > C_k(1-k^{-1/2})$,  then
\[(\sig_k^{(3)})^{-1} \leq (\frac{3}{2}+o(1)) k^2 \log k.\] 
(Although we will not apply this, we note that a clever iteration method allows for a slightly more effective choice of $\eta(s,k)$ and leads to very slightly better values for $\sig_k^{(2)}$ and $\sig_k^{(3)}$ for certain small $k$; the interested reader is referred to \cite{Ford}.) 

\subsection{The final result}
Using the value of $\sig_k(\be_0)$ provided by Lemma \ref{sig_lemma} to bound $S_{2^j}(\theta)$ in (\ref{minor_total}), it follows that the contribution of the minor arcs is in $L^{r,\infty}[0,1]$ (in fact in $L^\infty[0,1]$) if $\lam > 1-\sig_k(\be_0)$. Balancing this with the restriction (\ref{maj_assp}) arising in the major arcs that $\lam>\be_0$,  shows that $m_{k,\lam} \in L^{r,\infty}[0,1]$ for $\lam> \min (\lam_k^{(1)},\lam_k^{(2)},\lam_k^{(3)} ),$ where
\begin{eqnarray}
\lam_k^{(1)} & = & 1 - \frac{1}{2^{k-1}+1} \nonumber \\
\lam_k^{(2)} & = & 1 - \frac{1-(k-1)^{-2}}{4(k-1)^2 \log(k-1)+1} \label{lam_choices}\\
\lam_k^{(3)} & \approx & 1 - \frac{1}{(\frac{3}{2}+o(1))k^2 \log k} \nonumber,
\end{eqnarray}
where the value $\lam_k^{(3)}$ as stated is for appropriately large $k$.
(For comparison, the original result of Proposition 3 of \cite{SW2} stated that $m_{k,\lam} \in L^{r,\infty}[0,1]$ for $\lam$ ``sufficiently close'' to 1; although Stein and Wainger did not state the range of $\lam$ explicitly, in their setup $\be=k$, $\be_1=1$, and $\be_0=1/2$, so that applying the Weyl bound shows that $m_{k,\lam} \in L^{r,\infty}[0,1]$ in the smaller range $1-2^{-k} < \lam < 1$.)

To summarise, as a consequence of (\ref{lam_choices}) and Lemma \ref{folk_lemma}, we have proved:
\begin{thm}\label{thm_circle}
Let 
\[ \lam_k= \min (\lam_k^{(1)},\lam_k^{(2)},\lam_k^{(3)} ),\]
where the $\lam_k^{(t)}$ are as in (\ref{lam_choices}).
Then the operator $I_{k,\lam}$ extends to a bounded operator from $\ell^p(\Z)$ to $\ell^q(\Z)$ with
\[ ||I_{k,\lam} f||_{\ell^q} \leq A ||f||_{\ell^p} \]
 if $\lam_k < \lam < 1$ and $p,q$ satisfy
\newline
 (i) $1/q \leq 1/p - (1-\lam)/k$
\newline
 (ii)  $1/q<\lam, 1/p >1-\lam.$  
\end{thm}
Numerically, the first option $\lam_k^{(1)}$ is the minimum for $1\leq k \leq 11$.
For example, this theorem proves Conjecture \ref{Conj} for the operator $I_{3,\lam}$ in the range $4/5 < \lam < 1$.
The second option $\lam_k^{(2)}$ behaves asymptotically roughly like $1-\frac{1}{4k^2 \log k}$, and is  only the minimum for  $11\leq k \leq 13$. For $k \geq 14$, $\lam_k^{(3)}$ is the minimum, behaving approximately like $1-\frac{1}{\frac{3}{2}k^2 \log k}$ for large $k$.

\section{Relation to Waring's problem}\label{sec_Hyp_K}
Theorem \ref{thm_circle} proves Conjecture \ref{Conj} in a limited range $\lam_k < \lam < 1$. We now turn to a second approach that will allow us to prove results toward Conjecture \ref{Conj} in longer ranges of $\lam$. Let $r_{s,k}(l)$ denote the number of representations of a positive integer $l$ as a sum of $s$ summands, each a $k$th power of a positive integer: $l=n_1^k + \cdots + n_s^k$.  Waring's problem asks for the least $s$ such that all $l$ have such a representation, as well as for asymptotics for $r_{s,k}(l)$. This problem has  been studied extensively by means of the circle method, which only gives a meaningful asymptotic for $r_{s,k}(l)$ if the number $s$ of variables is sufficiently large in relation to the degree $k$. For $s=k$, the asymptotic behaviour of $r_{k,k}(l)$, or even its average behaviour, remains a very difficult open question.

For $k=2$, it is trivially true that $r_{2,2}(l) = O(l^\ep)$ for any $\ep>0$. (In fact the number of representations of a positive integer $l$ by any positive definite binary quadratic form is $O(l^\ep)$.) Hardy and Littlewood \cite{HL} conjectured (Hypothesis $K$) that $r_{k,k}(l)=O( l^\ep)$ for every $k \geq 2$, but this has been shown to be false for $k =3$ by Mahler \cite{Mah}, and is expected to be false for $k \geq 4$ as well. However, such a result is still expected to hold on average, leading to the following  conjecture (termed Hypothesis $K^*$ by Hooley \cite{Hoo}):

\begin{conjecture}[(Hypothesis $K^*$)]\label{Conj_K}
For $k =2$, 
\[\sum_{l=1}^N (r_{k,k}(l))^2= O(N^{1+\ep})\]
as $N \maps \infty$, for any $\ep>0$.
For $k \geq 3$,
\beq\label{K_sum}
\sum_{l=1}^N (r_{k,k}(l))^2= O(N)
\eeq
as $N \maps \infty$.
\end{conjecture}

For $k=2$, this result follows trivially from the pointwise bound for $r_{2,2}(l)$. The case $k=3$ already presents substantial difficulties, and has only been proved conditionally on certain standard conjectures and the Riemann Hypothesis for Hasse-Weil $L$-functions by Heath-Brown \cite{HB98} and Hooley \cite{Hoo}, independently. For $k \geq 4$, Hypothesis $K^*$ remains unproved even conditionally. The trivial bound for the sum in (\ref{K_sum}) is $O(N^2)$, and currently all known results toward Hypothesis $K^*$ are of the form $O(N^{2 - \del_k})$ for some small $\del_k>0$ tending to zero as $k \maps \infty.$

Stein and Wainger noted in \cite{SW2} that Hypothesis $K^*$ is equivalent to the statement that $m_{k,\lam} \in L^{2k}[0,1]$ for all $1/2< \lam < 1$. 
We now prove, by a similar argument, that a far more general equivalence exists between $m_{k,\lam}$ and  $r_{s,k}(l)$.
For integers $s,k \geq 2$, let Property $K_{s,k}^*(\beta)$ denote the property that 
\beq\label{sk_sum}
\sum_{l=1}^N (r_{s,k}(l))^2 =O(N^{\beta}).
\eeq
(Thus for $k \geq 3$, Hypothesis $K^*$ is Property $K_{k,k}^*(1)$.)
The following equivalence holds:
\begin{prop}\label{m_Ks_equiv}
Given $\be>0$, the property that $m_{k,\lam} \in L^{2s}[0,1]$ for every $\beta k/2s < \lam < 1$ is equivalent to Property $ K_{s,k}^*(\beta).$
\end{prop}
First assume that Property $K_{s,k}^*(\beta)$ holds.
Define
\[ S_y(\theta) = \sum_{n=1}^\infty e^{-\pi n^k(y+2i\theta)}.\]
Note that
\[ n^{-\lam} = \pi^{\lam/k} \Ga(\lam/k) \int_0^\infty e^{-\pi n^ky} y^{\lam/k-1}dy,\]
so that 
\[ m_{k,\lam}(\theta)= c_{k,\lam} \int_0^\infty S_y(\theta) y^{\lam/k-1}dy.\]
The contribution of $\int_1^\infty$ is $O(1)$, thus
\[ ||m_{k,\lam}||_{L^{2s}[0,1]} \leq c_{k,\lam} \int_0^1 ||S_y(\cdot)||_{L^{2s}[0,1]} y^{ \lam/k-1}dy + O(1).\]
By Parseval's identity,
\[ \int_0^1 |S_y(\theta)|^{2s} d\theta = \sum_{l=1}^\infty r_{s,k}^2(l) e^{-2\pi ly}.\]
Thus assuming Property $K_{s,k}^*(\beta)$ holds, $||S_y(\cdot)||_{L^{2s}} = O(y^{-\beta/2s})$ and hence $m_{k,\lam} \in L^{2s}[0,1]$ for all $\beta k/2s < \lam <1$.

Conversely,
if $m_{k,\lam} \in L^{2s}[0,1]$ for all $\beta k/2s < \lam <1$, then  
\[ \left( 1 + \sum_{n=1}^\infty \frac{e^{-2\pi i n^k \theta}}{n^\lam} \right) ^s = \sum_{l=0}^\infty a_l e^{-2\pi i l \theta} \in L^2 [0,1]\]
for every $\beta k/2s < \lambda< 1$, where now
\[ a_l = \sum_{j \leq s} \sum_{l = n_1^k + \cdots + n_j^k} \frac{1}{n_1^\lam \cdots n_j^\lam}.\]
Note that $n_i^k \leq l$ for every $i$, so $n_1^\lam \cdots n_j^\lam \leq l^{j \lam /k} \leq l^{s\lam /k}$ and thus picking out the $j=s$ term, $a_l \geq r_{s,k}(l)l^{-s\lam/k}$. Thus the fact that $\sum_l a_le^{-2\pi i l \theta} \in L^2[0,1]$ implies by Parseval's identity that $\sum_{l=1}^\infty r_{s,k}^2(l)l^{-2s\lam/k} < \infty $ for all $\lam > \beta k/2s$, giving Property $K_{s,k}^*(\beta)$.\xtra{
		For in fact, partial summation shows that in order for this infinite sum to be finite, we must have 
\[\sum_{l=1}^N r_{s,k}^2(l)  = O(N^\ga)\] 
where $\ga < 2s\lam/k$ for all $\beta k/2s < \lam < 1$; thus we can take $\ga \leq 2s(\beta k /2s)/k = \beta$ as desired.}

The equivalence in Proposition \ref{m_Ks_equiv} shifts the study of the multiplier $m_{k,\lam}$ to the study of Property $K^*_{s,k}(\beta)$. The flexibility of the equivalence is crucial: while little is known toward Hypothesis $K^*$, very deep results are known for mean values of the form (\ref{sk_sum}) for $s$ sufficiently small or sufficiently large relative to $k$. Results obtained for $m_{k,\lam}$ via this equivalence can be translated  into results for the operator $I_{k,\lam}$ via Lemma \ref{K_to_op}.
In particular, assuming the truth of Hypothesis $K^*$ would allow the choices $u=2k$ and $\eta=1/2$ in Lemma \ref{K_to_op}, leading to Conjecture \ref{Conj} but with a strict inequality in condition (i) and in the range $1/2< \lam< 1$. (As noted before, in the limit as $k \maps \infty$, this is nearly the full range $\lam_k^*<\lam <1$ required to prove Conjecture \ref{Conj}.)
In general, it is clear that bounds obtained for $I_{k,\lam}$ from a particular Property $K_{s,k}^*(\beta)$ will never attain equality in condition (i), and hence will always  be weaker than Conjecture \ref{Conj}. However, the advantage in pursuing such bounds is that they will hold for a longer range of $\lam$ than the original approach of Theorem \ref{thm_circle}.

\section{Mean values of exponential sums}\label{sec_mean_vals}
We now record the best known results for Property $K_{s,k}^*(\beta)$. Let 
\beq\label{s_dfn}
S_{k,X}(\al) = \sum_{n=1}^X e^{2 \pi i n^k \al}
\eeq
and set 
\[\mathbf{I}_{s,k}(X) = \int_0^1 |S_{k,X}(\al)|^{2s} d\al.\]
Note that
 \beq\label{sm_I}
  \sum_{l=1}^N r_{s,k}^2(l) = \mathbf{I}_{s,k}(N^{1/k}) ,
  \eeq
so that Property $K^*_{s,k}(\beta)$ amounts to a bound for the mean value $\mathbf{I}_{s,k}(X)$.
The expected behaviour of $\mathbf{I}_{s,k}(X)$ depends on the size of $s$ with respect to $k$, and so we will consider three cases: $s< k$, $s=k$, and $s >k$.

\subsection{Property $K_{s,k}^*(\beta)$ for $s$ small}
Consider $2 \leq s<k$. Note that $\mathbf{I}_{s,k}(X)$ may also be expressed as the number of $(2s)$-tuples $(x_1, \ldots, x_{2s})$ of positive integers $x_i \leq X$ satisfying the equation 
\beq\label{2eqn}
x_1^k + \cdots + x_s^k = x_{s+1}^k + \cdots + x_{2s}^k.
\eeq
A solution is considered diagonal if $(x_{s+1}, \ldots, x_{2s})$ is a permutation of $(x_1, \ldots, x_{s})$; the total number of diagonal solutions of (\ref{2eqn}) is asymptotically of size $s! X^s.$ The \emph{paucity conjecture} (see \cite{SalWoo}) posits that for $s < k$, these trivial diagonal solutions dominate the solutions to (\ref{2eqn}),  so that for $2 \leq s<k$, it is conjectured that
\beq\label{N_bound}
\mathbf{I}_{s,k}(X) = s! X^s + o(X^s).
\eeq
This would imply Property $K_{s,k}^*(s/k)$ for $2 \leq s < k$.

For $s=2,$ $k >2,$ it may be proved quite simply that $\mathbf{I}_{s,k}(X) = O(X^{2+\ep})$, since the equation (\ref{2eqn}) takes the form 
\[ x_1^k - y_1^k = y_2^k - x_2^k\]
 and $I_{s,k}(X)$ may be regarded as counting simultaneous solutions to $l=x^k - y^k = (x-y)(x^{k-1} + \cdots + y^{k-1})$, of which there are very few since there are very few divisors of $l$. (This is proved rigourously, for example, in Proposition 4 of \cite{SW2}, which proves as a consequence that $m_{k,\lam} \in L^4[0,1]$ for all $1/2 < \lam < 1$ and $k \geq 2$.) For $s \geq 3$ the problem is substantially harder. In the specific case $s=3$, Browning and Heath-Brown \cite{BHB04} showed that (\ref{N_bound}) holds for $k>32$; this was later improved (also for $s=3$) by Salberger \cite{Sal05} to $k>25$. Recently, Salberger and Wooley \cite{SalWoo} have made substantial progress toward the paucity conjecture, showing that for any $k \geq 3$, the diagonal solutions dominate and (\ref{N_bound}) holds for $s \leq (\frac{1}{4} + o(1)) \log k ( \log \log k)^{-1}.$ As a consequence, Proposition \ref{m_Ks_equiv} and Lemma \ref{K_to_op} imply:

\begin{thm}\label{I_sk_small}
The operator $I_{k,\lam}$ extends to a bounded operator from $\ell^p (\Z)$ to $\ell^q(\Z)$ with 
\[ ||I_{k,\lam} f||_{\ell^q} \leq A ||f||_{\ell^p} \]
if  $1/2 < \lam < 1$ and  $p,q$ satisfy the conditions
\newline
(i-a) $1/q < 1/p - \frac{(1-\lam)}{s_k}$,
\newline
 (ii-a) $1/q < \lam, 1/p > 1-\lam$,\\
 where 
 \[ s_k \leq (\frac{1}{4} + o(1)) \log k ( \log \log k)^{-1}.\]
\end{thm}
This result holds in nearly all the desired range of $\lam$ (namely $\lam_k^*< \lam < 1$), but condition (i-a), roughly of the form 
\[ 1/q < 1/p - \frac{4 \log \log k(1-\lam)}{\log k}\]
for sufficiently large $k$,
 is significantly weaker than condition (i) of Conjecture \ref{Conj}.

\subsection{Property $K_{k,k}^*(\beta)$}\label{sec_kk}
Recall from Conjecture \ref{Conj_K} that for $s=k=2$, Property $K_{k,k}^*(1+\ep)$ known to be true, and for $s=k\geq 3$, Property $K_{k,k}^*(1)$ is expected to be true. For $k \geq 3$, the current best known results toward Property $K_{k,k}^*(\theta)$ follow from the classical lemma of Hua (see for example \cite{Hua}, \cite{IK}), which states that 
\[ \mathbf{I}_{2^{l-1},k}(X) =O( X^{2^l - l + \ep}),\]
for any integer $1 \leq l \leq k$, with implied constant depending on $k,\ep$.
Thus if $k$ is a power of $2$, using the notation (\ref{s_dfn}),
\[ \sum_{l=1}^N r_{k,k}^2(l)= \int_0^1 |S_{k,N^{1/k}}(\al)|^{2k}  d\al= O(N^{2 - ( \log_2 k+1)/k + \ep}),\]
which gives Property $K_{k,k}^*(2 - \del_k)$ with $\del_k = (\log_2 k+1)/k$.
Incidentally, since trivially $|S_{k,X}(\al)| \leq X$,  this argument also shows that Property $K_{k,k}^*(2)$ is trivial for all $k \geq 2$. 

If $k$ is not a power of $2$, one can use H\"{o}lder's inequality to interpolate between two applications of Hua's lemma. For example, in the case $k=3$, 
\[ \int_0^1 |S_{3,N^{1/3}}(\al)|^6 d\al \leq \left(\int|S_{3,N^{1/3}}|^4\right)^{1/4} \left( \int |S_{3,N^{1/3}}|^8 \right)^{1/8}.\]
Applying Hua's inequality to each term on the right gives $\int |S_{3,N^{1/3}}|^4 \ll X^2$  and $\int |S_{3,N^{1/3}}|^8 \ll X^5$  and hence 
\[\sum_{l=1}^N r_{3,3}^2(l)  = O(N^{7/6 + \ep}).\]
More generally, the same principle yields:
\begin{prop}\label{Hua-Holder}
If $2^{s} \leq k < 2^{s+1},$ then Property $K_{k,k}^*(2 -\del_k)$ holds with $\del_k= (k + 2^{s}s)/(2^{s}k)$.
\end{prop}
This yields the best known unconditional bounds toward Hypothesis $K^*$ for $k \geq 3$; for example, this proves Property $K_{3,3}^*(7/6)$ and $K_{4,4}^*(5/4)$.
Taking $u = 2k$ and $\eta = 1-\del_k/2$ in Proposition \ref{m_Ks_equiv} and Lemma \ref{K_to_op} proves:

\begin{thm}\label{I_kk}
For $2^s \leq k <2^{s+1}$, let
\[\lam_k =1 -   \frac{k + 2^{s}s}{2^{s+1}k}.\]
Then the operator $I_{k,\lam}$ extends to a bounded operator from $\ell^p (\Z)$ to $\ell^q(\Z)$ with 
\[ ||I_{k,\lam} f||_{\ell^q} \leq A ||f||_{\ell^p} \]
if  $\lam_k< \lam < 1$ and  $p,q$ satisfy the conditions
\newline
(i-b) $1/q < 1/p - \frac{(1-\lam)}{2k(1-\lam_k)}$,
\newline
 (ii-b) $1/q < \lam, 1/p > 1-\lam$. 
\end{thm}

For example, this shows that $I_{3,\lam} : \ell^p \maps \ell^q$ for $1/q < 1/p - \frac{2(1-\lam)}{5}$ in the range  $7/12 < \lam < 1$. In general, for $k$ a power of $2$, $I_{k,\lam}$ maps $\ell^p \maps \ell^q$ for $1/q < 1/p - \frac{(1-\lam)}{\log_2 k +1}$ in the range $1 - \frac{\log_2 k+1}{2k}< \lam < 1$. This is intermediate between Theorem 1 and Theorem 2, both in terms of the range of allowable $\lam$ and in that condition (i-b) is stronger than (i-a) but weaker than (i).

\subsection{Property $K^*_{s,k}(\beta)$ for $s$ large}
 If $s \geq k \geq 3$, simple estimates (see \cite{Ford}) show that 
\beq\label{I_lowerbd}
\mathbf{I}_{s,k}(X) \gg X^{2s-k},
\eeq
and this is conjectured to be the true order of magnitude of $\mathbf{I}_{s,k}(X)$ in this range. If true, this would imply Property $K_{s,k}^*(2s/k-1)$. 
If $s$ is sufficiently large with respect to $k$, the circumstances are more favourable than for $s$ close to $k$. 
A clever argument of Ford \cite{Ford} uses existing bounds for ${\bf J}_{s,k}(N)$ (as defined in (\ref{J_eta})) to prove
\beq\label{Ford_bd}
 \mathbf{I}_{s,k}(X) \ll X^{2s - k + (\sqrt{2e}/k) \eta(s,k)},
 \eeq
 and hence Property $K_{s,k}^*(2s/k-1 + \del_k)$, where $\del_k = k^{-1}(\sqrt{2e}/k) \eta(s,k).$
We recall from (\ref{eta_sk}) that  the best known value is
$\eta(s,k)\approx k^2 e^{-2s/k^2}$, so that (\ref{Ford_bd}) is very close to the desired order of magnitude $X^{2s-k}$ as soon as $s$ is of the order $k^2 \log k$. 

Applying this result to the operator $I_{k,\lam}$ leads to a family of $\ell^p \maps \ell^q$ bounds intermediate between Theorems \ref{thm_circle} and \ref{I_kk}, with the precise results depending on the relative sizes of $s,k$. 
However, for clarity of presentation and maximal contrast to Theorem \ref{I_kk}, we will only state the theorem for $I_{k,\lam}$ resulting from applying (\ref{Ford_bd}) to derive asymptotics of the correct order of magnitude for $r_{s,k}(l)$.

 The circle method proves that for $s$ sufficiently large with respect to $k$,
\beq\label{r_asymp}
 r_{s,k}(l) = \frac{\Gamma(1+1/k)^s  }{\Gamma(s/k)}\mathfrak{S}(l)l^{s/k-1}+ O(l^{s/k-1-\del})
 \eeq
for an explicit constant $\del>0$, where $\mathfrak{S}(N)$ is the singular series and the implied constant in the error term depends only on $s,k$. (For precise definitions, see Theorem 20.2 of \cite{IK}.) This is only a true asymptotic if the singular series does not vanish, but as we will apply this result in the form of an upper bound for $r_{s,k}(l)$, we will disregard the possible vanishing of the singular series.

Let $\tilde{G}(k)$ denote the smallest $s$ for which the asymptotic (\ref{r_asymp}) holds. Then 
for $s \geq \tilde{G}(k)$, Property $K_{s,k}^*(2s/k-1)$ holds, from which our ultimate bound for the operator $I_{k,\lam}$ will follow.
The strength of Ford's bound (\ref{Ford_bd}) in this context is that it pushes the allowable range of $s$ closer to $k$. 
It is conjectured that it should be possible to take $\tilde{G}(k)$ on the order of $k$ (leading back toward Hypothesis $K^*$). The original result of Hardy and Littlewood \cite{HL22} in 1922 stated that $\tilde{G}(k) \leq (k-2) 2^{k-1} +5$. Since then, $\tilde{G}(k)$ has been whittled down gradually; for large $k$ key historical results are due to Hua \cite{Hua38} \cite{Hua49}, Vinogradov \cite{Vin85}, Wooley \cite{Woo92} and Ford \cite{Ford}, while for small $k$, key results are due to Kloosterman \cite{Kloos}, Vaughan \cite{Vau86A} \cite{Vau86B}, Heath-Brown \cite{HB88}, and Boklan \cite{Bok}.

For $3 \leq k\leq 8$, the current best known results are as follows: $\tilde{G}(3) \leq 8$,  $\tilde{G}(4) \leq 16$, $\tilde{G}(5) \leq 32$ (Vaughan \cite{Vau86A}, \cite{Vau86B}),  $\tilde{G}(6) \leq 56$, $\tilde{G}(7) \leq 112$, $\tilde{G}(8) \leq 224$ (Boklan \cite{Bok}). For $k \geq 9$,
Ford's work \cite{Ford} based on the bound (\ref{Ford_bd}) currently holds the record: for example, $\tilde{G}(9) \leq 393$, and asymptotically $\tilde{G}(k)  \leq k^2(\log k + \log\log k +O(1))$.
Applying these results to the operator $I_{k,\lam}$ by taking $u = 2s$ and $\eta = 1-k/2s$  in Lemma \ref{K_to_op} yields our final result:
\xtra{We also obtain a partial result for $0<\lam \leq \nu_k$ by interpolating the result of Proposition \ref{K_to_op} with  the trivial $\ell^1 \maps \ell^\infty$ bound for $I_{k,\lam}$ whenever $\Re(\lam)\geq 0$; in this respect we note that $ -1 + \lam(1-1/r)/\eta $ is the linear function in $\lam$ which gives $-1/r $ when $\lam = \eta $ and $-1$ when $\lam=0$. 
\begin{cor}
For any $s \geq \tilde{G}(k)$, $I_{k,\lam} : \ell^p \maps \ell^q$ for all $0 < \lam \leq 1-k/2s$ with  $p,q$ satisfying the conditions
\newline
(i) $1/q < 1/p  -1 +  \frac{\lam(1-1/2s)}{(1-k/2s)} $
\newline
(ii)$1/q < \lam, 1/p > 1-\lam$.
\end{cor}
\ques{Why does this appear to get stronger as $s$ is bigger? The range would also change.}
}
\begin{thm}\label{I_sk}
Let 
\[\lam_k = 1-\frac{k}{2\tilde{G}(k)},\]
 where $\tilde{G}(k)$ is the least $s$ for which the asymptotic (\ref{r_asymp}) holds. 
Then the operator $I_{k,\lam}$ extends to a bounded operator from $\ell^p(\Z)$ to $\ell^q(\Z)$ with 
\[ ||I_{k,\lam}f||_{\ell^q} \leq A ||f||_{\ell^p} \]
 if $\lam_k < \lam < 1$ and $p,q$ satisfy the conditions
\newline
(i-c) $1/q < 1/p - \frac{(1-\lam)}{k}$,
\newline
 (ii-c) $1/q < \lam, 1/p > 1-\lam$.\\

\end{thm}
This theorem illustrates the trade-off inherent in this approach: requiring $s$ to be sufficiently large that sharp asymptotics hold for $r_{s,k}(l)$ limits the allowable range of $\lam$ for the operator bounds.
 In particular, the range of allowable $\lam$ in Theorem \ref{I_sk} is shorter than that of either Theorem \ref{I_sk_small} or \ref{I_kk}, although it is longer than the range allowed by Theorem \ref{thm_circle}. For example, Theorem \ref{I_sk} shows that $I_{3,\lam}: \ell^p \maps \ell^q$ for $1/q< 1/p - (1-\lam)/3$ in the range $13/16 < \lam < 1$; in the limit, the allowable range is approximately $1 - \frac{1}{2k \log k} < \lam < 1$.  
On the other hand, condition (i-c) is stronger than both conditions (i-a) and (i-b) of Theorems \ref{I_sk_small} and \ref{I_kk}.  Indeed condition (i-c) is nearly  (within $\ep$) sharp, compared to the optimal condition (i) of Conjecture \ref{Conj}. 

In conclusion, direct approaches to the multiplier $m_{k,\lam}$ via the circle method depend on the best known bounds for individual Weyl sums, and lead to results for $I_{k,\lam}$ that are sharp with respect to the relationships between $p,q$ and $\lam$, but hold for a very small range of $\lam$. Indirect approaches to $m_{k,\lam}$ via mean values of Weyl sums greatly extend the allowable range of $\lam$, but provide weaker relations between $p,q$ and $\lam$. This phenomenon is visible in Figure \ref{fig_1}.  It is worth noting that any improvements to bounds for $\mathbf{J}_{s,k}$ and $\mathbf{I}_{s,k}$ will have immediate impact on the operator $I_{k,\lam}$.

\bibliographystyle{amsplain}
\bibliography{AnalysisBibliography}


\end{document}

%% file: format.tex
\setcounter{secnumdepth}{4}


\newtheorem{thm}{Theorem}
\newtheorem{prop}{Proposition}
\newtheorem{lemma}[prop]{Lemma}
\newtheorem{cor}{Corollary}[thm]
\newtheorem{conjecture}{Conjecture}

\newtheorem{letterthm}{Theorem}
\renewcommand{\theletterthm}{\Alph{letterthm}}

\newtheorem{letterlemma}{Lemma}
\renewcommand{\theletterlemma}{\Alph{letterlemma}}

%

\renewcommand{\labelenumi}{(\roman{enumi})}

\newcommand{\Leg}[2]{ \left( \frac{#1}{#2} \right)}
\newcommand{\con}{\equiv}
\newcommand{\ncon}{\not\equiv}
\newcommand{\ndiv}{\nmid}
\newcommand{\modd}[1]{\; ( \text{mod} \; #1)}
\newcommand{\bstack}[2]{\substack{#1 \\ #2}}
\newcommand{\tstack}[3]{#1\atop {#2 \atop #3}}
\newcommand{\sprod}[2]{\langle #1, #2 \rangle}
\newcommand{\inv}[1]{\overline{#1}}
\newcommand{\conj}[1]{\overline{#1}}
\newcommand{\maps}{\rightarrow}
\newcommand{\lequiv }{\Longleftrightarrow}
\newcommand{\intersect}{\cap}
\newcommand{\Intersect}{\bigcap}
\newcommand{\bigunion}{\bigcup}
\newcommand{\Union}{\bigcup}
\newcommand{\union}{\cup}
\newcommand{\grad}{\nabla}
\newcommand{\diag}{\text{diag}}
\newcommand{\uhat}[1]{\check{#1}}
\newcommand{\tensor}{\otimes}
\newcommand{\isom}{\cong}
\newcommand{\deriv}[1]{\frac{\partial}{\partial #1}}
\newcommand{\pderiv}[2]{\frac{\partial ^#1}{\partial #2}}
\newcommand{\comp}{\circ}

\newcommand{\Frac}{\text{Frac}}
\newcommand{\Ker}{\text{Ker}}
\newcommand{\poly}{\text{poly}}
\newcommand{\Poly}{\text{Poly}}
\newcommand{\Id}{\text{Id}}
\newcommand{\lcm}{\text{lcm}}

\newcommand{\al}{\alpha}
\newcommand{\be}{\beta}
\newcommand{\gam}{\gamma}
\newcommand{\ga}{\gamma}
\newcommand{\del}{\delta}
\newcommand{\ep}{\epsilon}
\newcommand{\ka}{\kappa}
\newcommand{\om}{\omega}
\newcommand{\Om}{\Omega}
\newcommand{\Sig}{\Sigma}
\newcommand{\sig}{\sigma}
\newcommand{\lam}{\lambda}
\newcommand{\Lam}{\Lambda}
\newcommand{\Ga}{\Gamma}

\newcommand{\uchi}{\underline{\chi}}
\newcommand{\tchi}{\tilde{\chi}}

\newcommand{\A}{\mathcal{A}}
\newcommand{\Acal}{\mathcal{A}}
\newcommand{\Bcal}{\mathcal{B}}
\newcommand{\Dcal}{\mathcal{D}}
\newcommand{\E}{\mathcal{E}}
\newcommand{\Ecal}{\mathcal{E}}
\newcommand{\Fcal}{\mathcal{F}}
\newcommand{\Hcal}{\mathcal{H}}
\newcommand{\I}{\mathcal{I}}
\newcommand{\Ical}{\mathcal{I}}
\newcommand{\J}{\mathcal{J}}
\newcommand{\Jcal}{\mathcal{J}}
\newcommand{\K}{\mathcal{K}}
\newcommand{\Lcal}{\mathcal{L}}
\newcommand{\Mcal}{\mathcal{M}}
\newcommand{\Ncal}{\mathcal{N}}
\newcommand{\Rcal}{\mathcal{R}}
\newcommand{\Scal}{\mathcal{S}}
\newcommand{\Tcal}{\mathcal{T}}
\newcommand{\U}{\mathcal{U}}
\newcommand{\Ut}{\scriptscriptstyle{\mathcal{U}}}
\newcommand{\Vcal}{\mathcal{V}}
\newcommand{\V}{\mathcal{V}}

\newcommand{\C}{\mathbb{C}}
\newcommand{\F}{\mathbb{F}}
\newcommand{\HH}{\mathbb{H}}
\newcommand{\N}{\mathbb{N}}
\newcommand{\Proj}{\mathbb{P}}
\newcommand{\Q}{\mathbb{Q}}
\newcommand{\R}{\mathbb{R}}
\newcommand{\T}{\mathbb{T}}
\newcommand{\Z}{\mathbb{Z}}

\newcommand{\p}{\mathfrak{p}}
\newcommand{\Nf}{\mathfrak{N}}
\newcommand{\M}{\mathfrak{M}}
\newcommand{\m}{\mathfrak{m}}

\newcommand{\x}{{\bf x}}
\newcommand{\y}{{\bf y}}

\newcommand{\res}{\text{res}}
\newcommand{\Tr}{\text{Tr}}
\newcommand{\rank}{\text{rank}\,}
\newcommand{\cond}{\text{cond}\,}
\newcommand{\Disc}{\text{Disc}}

\newcommand{\beq}{\begin{equation}}
\newcommand{\eeq}{\end{equation}}